\documentclass{amsart}

\usepackage[latin1]{inputenc}
\usepackage{amsmath}
\usepackage{amsfonts}
\usepackage{amssymb}
\usepackage[square,numbers]{natbib}
\usepackage{enumerate}
\usepackage[small]{diagrams}

\usepackage{bbold}

\newcommand\id{\mathbb 1}
\newcommand\Id{\operatorname{Id}}
\newcommand\Hom{\operatorname{Hom}}
\newcommand\St{\operatorname{St}}
\newcommand\End{\operatorname{End}}
\newcommand\Ext{\operatorname{Ext}}

\newtheorem{thm}{Theorem}

\newtheorem{alg}[thm]{Algorithm}

\begin{document}

%\begin{frontmatter}
  \title{Blackbox computation of $A_\infty$-algebras}
  \author{Mikael Vejdemo-Johansson}
  \thanks{Partially supported by DFG Sachbehilfe grant GR 1585/4-1}
  %\ead{mik@math.uni-jena.de}
\email{mik@math.uni-jena.de}
  \address{Friedrich-Schiller-Universität Jena, Germany\\Currently at
    Stanford University}
  %\corauth{Mikael Vejdemo-Johansson}
    
\dedicatory{To Tornike Kadeishvili}

\maketitle

   \begin{abstract}
    Kadeishvili's proof of the minimality theorem (Kadeishvili, T.;
    \textit{On the homology theory of fiber spaces}, Russian
    Math. Surveys 35 (3) 1980)
    induces an algorithm for the inductive computation of an
    $A_\infty$-algebra structure on the homology of a dg-algebra.

    In this paper, we prove that for one class of dg-algebras, the
    resulting computation will generate a complete $A_\infty$-algebra
    structure after a finite amount of computational work.
  \end{abstract}

%\begin{keyword}

  %\end{keyword}
%\end{frontmatter}

\section{Introduction}
\label{sec:introduction}

$A_\infty$-algebras have been in use in topology since their
introduction by Stasheff \citep{stasheff_ha_hspaces63}. Their
applicability in an algebraic context was made clear by the Minimality
Theorem, proven by Kadeishvili \citep{kadeishvili80}, then extended and
reproven using techniques both from homological perturbation theory
(see \citep{johansson-lambe-01} for a good overview and
\citep{huebschmann_homotopy_1986} for one origin of the techniques) and explicit
tree formulae \citep{huebschmann-kadeishvili-91,merkulov_sha_kahler99}.

These structures are found in representation theory as well:
Using suitable $A_\infty$-structures, Huebschmann has computed various
group cohomology rings of groups given as extensions; the spectral
sequence of the group extension under consideration is then an
invariant of the $A_\infty$-structure \citep{huebschmann-coho-89,huebschmann-modp-89,huebschmann-perturbation-89,huebschmann-91}.
Furthermore both Keller \citep{keller_intro01,keller_addendum_2001,keller_Aoo_repr00} and Lu,
Palmieri, Wu and Zhang
\citep{lu_wu_palmieri_zhang_Aoo_ring04,lu_wu_palmieri_zhang_Aoo_Ext06}
have studied the use of $A_\infty$-algebras for module categories and
cohomology rings over specific algebras.

In many of the proofs of the Minimality theorem, a more or less
explicitly stated version of a homological perturbation theory
approached were used. This version uses decorated trees to give
explicit formulas for the higher operations in terms of a deformation
retract and the multiplication map from the original dg-algebra. The
technique is present in the work of Huebschmann and Kadeishvili, and
has very explicit treatments by Merkulov and by Lu, Palmieri, Wu and
Zhang \citep{huebschmann-kadeishvili-91,merkulov_sha_kahler99,
  lu_wu_palmieri_zhang_Aoo_Ext06,lu_wu_palmieri_zhang_Aoo_ring04}. More
details on the history of $A_\infty$-algebras may be found in
Huebschmann's survey talk \citep{huebschmann-origins-07} as well as in
\citep{huebschmann_construction_2008}.

The techniques used to prove the Minimality theorem all yield explicit
methods to compute an $A_\infty$-algebra structure, and the main known
techniques for producing such structures are rooted in the various
proofs available. All the methods rooted in homological perturbation
theory yield complete structures, but at a price: they all require
global information of a kind that in algebraic settings is not
necessarily easy to come by. By contrast, the Kadeishvili algorithm,
derived directly from Kadeishvili's original proof, only gives
fragments of the $A_\infty$-structure, but does so without the kind of
global information expected in the other methods. We merely need some
sort of \emph{blackbox}, capable of performing computations within the
dg-algebra we start with to extract information about the resulting
$A_\infty$-algebra structure on its homology.

In this paper, we prove the following theorems, which together yield a
strong reduction in the computational load in computing
$A_\infty$-algebra structures on particularly nice rings:

\newtheorem*{thma}{Theorem \ref{lemma:periodic}}
\newtheorem*{thmb}{Theorem \ref{lemma:halts}}

\newcommand{\theorema}[1]{
  Suppose that $A$ is a dg-algebra and suppose that
  \begin{enumerate}[{A}1.]
  \item there is an element $z\in H_*A$ generating a
    polynomial subalgebra and that $H_*A$ is free as a
    $\mathbb{k}[z]$-module and
  \item\label{#1:2} $H_*A$ possesses a $\mathbb{k}[z]$-linear $A_{n-1}$-algebra
	structure induced by the dg-algebra structure on $A$, such that
    \[
	f_1(z)f_k(a_1,\dots,a_k)=f_k(a_1,\dots,a_k)f_1(z)
    \]
    
  \item $b_1,\dots$ is a $\mathbb{k}[z]$-basis of $H_*A$, and that we have
    chosen $m_n(v_1,\dots,v_n)$ and $f_n(v_1,\dots,v_n)$ according to
    the Kadeishvili algorithm for all $v_i\in\{b_1,\dots\}$.
  \end{enumerate}

  Then a choice of $m_n(v_1,\dots,v_n)$ and $f_n(v_1,\dots,v_n)$
  according to the Kadeishvili algorithm for all $v_i$ taken from a
  $\mathbb{k}[z]$-basis for $H_*A$ extends to a $\mathbb{k}[z]$-linear
  $A_n$-algebra structure for $H_*A$ induced by the dg-algebra
  structure on $A$.
}
\begin{thma}
\theorema{thm1}
\end{thma}

%%% TODO: Hack to get numbering right.
\newcommand\theoremb{
  Let $A$ be a dg-algebra. Suppose that in the computation of an
  $A_{2q-2}$-structure on $H_*A$, we have shown that 
  $f_k=0$ and $m_k=0$ for all $q\leq k\leq 2q-2$ for some $q$.
  
  Then this $A_{2q-2}$-algebra structure extends to an
  $A_\infty$-algebra structure with $f_k=0$ and $m_k=0$ for all $k\geq
  q$.
}
\begin{thmb}
\theoremb
\end{thmb}

In the case where $H_*A$ is a finite $\mathbb{k}[z]$-module, this reduces the
workload necessary in order to inductively compute an $A_n$-algebra
structure on $H_*A$ to a finite workload. Additionally, we adapt the
conditions of Theorem \ref{lemma:periodic} in Theorem
\ref{lemma:commutes} to recompute the $A_\infty$-algebra structure on
the Yoneda Ext-algebra over the truncated polynomial ring originally
computed by Madsen \citep{madsen_phd}.

The paper is organized into five sections. Section
\ref{sec:introduction}, is this introduction. Section
\ref{sec:infinity-algebras} defines $A_\infty$-algebras and the
properties used. Section \ref{sec:kade-algor} introduces
the Minimality Theorem, Kadeishvili's proof thereof and the algorithm
inherent in the proof. In Section \ref{sec:comp-reduct}, the theorems
discussed above are stated and proven, and finally in Section
\ref{sec:implementation}, we use the Computer Algebra
System \textsc{Magma} to implement the techniques we introduced.

\section{A-infinity algebras}
\label{sec:infinity-algebras}

Let $\mathbb{k}$ be a field. An $A_\infty$-algebra is a graded
$\mathbb{k}$-vector space $A$ equipped with a family
$\{m_n\}_{n=1}^\infty$ of multilinear maps $m_n\in\Hom(A^{\otimes
n},A)$ of degree $2-n$ that satisfy the Stasheff family of identities
\[
\St_n: \sum_{r+s+t=n} (-1)^{r+st} m_{r+1+t}(\id^r\otimes m_s\otimes \id^t)=0
\]

A morphism $f:A\to B$ of $A_\infty$-algebras is a family
$\{f_n\}_{n=1}^\infty$ of multilinear maps $f_n\in\Hom(A^{\otimes
n},A)$ of degree $1-n$ that satisfy the Stasheff family of morphism
identities
\[
\St^m_n: 
\sum_{r+s+t=n}(-1)^{r+st}f_{r+1+t}(\id^r\otimes m_s\otimes\id^t) =
\sum(-1)^wm_q(f_{i_1}\otimes\dots\otimes f_{i_q})
\]
where $i_1+\dots+i_q=n$ and
$w=(q-1)(i_1-1)+\dots+2(i_{q-2}-1)+(i_{q-1}-1)$.

We also define an $A_n$-algebra to be a graded vector space with a
family $\{m_k\}_{k=1}^n$ of multilinear maps satisfying the Stasheff
identities $\St_1,\dots,\St_n$. A morphism of $A_n$-algebras is a
family $\{f_k\}_{k=1}^n$ of multilinear maps satisfying the Stasheff
identities $\St^m_1,\dots,\St^m_n$.

An $A_n$-, or $A_\infty$-algebra is a dg-algebra if $m_n=0$ for all
$n\geq 3$. 

For a dg-subalgebra $R\subseteq A$ of an $A_n$-, or
$A_\infty$-algebra, we call $A$ an $R$-linear $A_n$-, or
$A_\infty$-algebra if $m_k\in\Hom_R(A^{\otimes k},A)$ for all $k$. A
morphism $f_*:A\to B$ of $A_n$-, or $A_\infty$-algebras is called
$R$-linear if $f_1(R)$ is a dg-subalgebra of $B$ and for any $r\in
R$ and $a_1,\dots,a_k\in A$, 
\[
  f_k(a_1,\dots,ra_i,\dots,a_k) = f_1(r)f_k(a_1,\dots,a_i,\dots,a_k)
\]
or in other words, viewing $B$ as an $R$-module through the map $f_1$,
that $f_k\in\Hom_R(A^{\otimes k},B)$ for all $k$.

\section{Kadeishvili's algorithm}
\label{sec:kade-algor}

At the core of our approach to the computation of $A_\infty$-algebras
is the minimality theorem:

\begin{thm}[Minimality]
  Let $A$ be a dg-algebra so that $H_*A$ is a free $\mathbb k$
  module. There is an $A_\infty$-algebra structure on $H_*A$ and an
  $A_\infty$-algebra quasi-isomorphism $f:H_*A\to A$ such that $f_1$
  is a cycle-choosing quasi-isomorphism of dg-vector spaces, $m_1=0$
  and $m_2$ is induced by the multiplication in $A$.

  The resulting structure is unique up to isomorphism.

  If $A$ has a unit $1$, then the structure and quasi-isomorphism can be
  chosen to be strictly unital, i.e. such that 
  $\forall k,\lambda\in\mathbb k, a_i\in A$
  \[
  f_k(a_1,\dots,a_{i-1},\lambda\cdot1,a_{i+1},a_k)=0 \qquad
  m_k(a_1,\dots,a_{i-1},\lambda\cdot1,a_{i+1},a_k)=0 
  \]
\end{thm}
\begin{proof}
  The following proof is from Kadeishvili's original paper
  \citep{kadeishvili80}. We review most of it here as a basis for the
  algorithm.

  Since $A$ is a dg-algebra, there are only two operations in the
  $A_\infty$-structure on $A$: $m_1$ and $m_2$. We shall denote $m_1$
  by $d$ and $m_2$ by the infix operator $\cdot$ or by juxtaposition.

  To initialize an inductive definition of an $A_\infty$-structure on
  $H_*A$, we pick $m_1=0$ and $m_2$ the induced multiplication on the
  coclasses. Furthermore, we let $f_1$ be some cycle-choosing $\mathbb
  k$ module homomorphism $H_*A\to A$.

  Set $\Psi_2(a_1,a_2)=f_1(a_1a_2)-f_1(a_1)f_1(a_2)$. This is a
  boundary, since $f_1(a_1a_2)$ is defined to be a representative
  cycle of the homology class containing $f_1(a_1)f_1(a_2)$. Hence,
  there is some $w$ such that $dw=\Psi_2(a_1,a_2)$. We define
  $f_2(a_1,a_2)=w$.

  Now, for $n>2$, write 
  \begin{multline*}
    \Psi_n(a_1,\dots,a_n)=
    \sum_{s=1}^{n-1}(-1)^{\varepsilon_1(a_1,\dots,a_n,s)}
    f_s(a_1,\dots,a_s)\cdot f_{n-s}(a_{s+1},\dots,a_n) + \\
    \sum_{j=2}^{n-1}\sum_{k=0}^{n-j}(-1)^{\varepsilon_2(a_1,\dots,a_n,k,j)}
    f_{n-j+1}(a_1,\dots,a_k,m_j(a_{k+1},\dots,a_{k+j}),\dots,a_n)
  \end{multline*}
  where the expressions
  $\varepsilon_1(a_1,\dots,a_n,s)=s+(n-s+1)(|a_1|+\dots+|a_s|)$ and
  $\varepsilon_2(a_1,\dots,a_n,k,j)=k+j(n-k-j+|a_1|+\dots+|a_k|)$ are
  the signs in the Stasheff morphism axiom $\St^m_n$ with the Koszul
  signs introduced.
  
  This $\Psi_n$ is the complete expression of the Stasheff morphism axiom
  $\St^m_n$, but with the two terms $f_1m_n$ and $m_1f_n$ removed. The
  central point of this proof is to fill in these terms. 

  By some tedious technical checking, we can confirm that the element
  $\Psi_n(a_1,\dots,a_n)\in\ker d$. Hence, $\Psi_n(a_1,\dots,a_n)$
  belongs to some coclass $z\in H_*A$. We define 
  \[
    m_n(a_1,\dots,a_n)=z\quad.  
  \]
  Since now $f_1(m_n(a_1,\dots,a_n))$ and $\Psi_n(a_1,\dots,a_n)$
  are in the same coclass, there is some coboundary $dw$, $w\in A$,
  such that $f_1(m_n(a_1,\dots,a_n))-\Psi_n(a_1,\dots,a_n)=dw$. We set
  \[
  f_n(a_1,\dots,a_n)=w\quad.
  \]
  Since we defined everything precisely in order to match the Stasheff
  axioms, we obtain with a structure that satisfies the Stasheff axioms.

  For example, we note that:
  \begin{align*}
  \Psi_3(a,b,c)&=(-1)^{\varepsilon_1(a,b,c,1)}f_1(a)f_2(b,c) + 
  (-1)^{\varepsilon_1(a,b,c,2)}f_2(a,b)f_1(c) + \\&\phantom{=}
  (-1)^{\varepsilon_2(a,b,c,0,2)}f_2(m_2(a,b),c) +
  (-1)^{\varepsilon_2(a,b,c,1,2)}f_2(a,m_2(b,c)) \\
  &=(-1)^{1+1\cdot|a|}f_1(a)f_2(b,c) + 
  (-1)^{2+2\cdot(|a|+|b|)}f_2(a,b)f_1(c) + \\&\phantom{=}
  (-1)^{0+2\cdot(\cdots)}f_2(m_2(a,b),c) +
  (-1)^{1+2\cdot(\cdots)}f_2(a,m_2(b,c)) \\
  =-(-1)^{|a|}&f_1(a)f_2(b,c) + f_2(a,b)f_1(c) + 
  f_2(m_2(a,b),c) - f_2(a,m_2(b,c))
  \end{align*}
  
  Regarding unitality, first we consider 
  $m_2(1,a) = m_2(a,1) = a$; then
  $\Psi_2(1,a)=a-a=0$ and $\Psi_2(a,1)=a-a=0$. Thus, we can safely
  choose $f_2(1,a)=f_2(a,1)=0$.

  Now, consider $\Psi_3$. We have three cases to consider:
  \begin{align*}
    &\Psi_3(1,a,b) \\&= 
    -(-1)^{|1|}f_1(1)f_2(a,b) + f_2(1,a)f_1(b) + 
    f_2(m_2(1,a),b) - f_2(1,m_2(a,b)) \\
    &= -f_2(a,b)+0+f_2(a,b)-0 = 0 \\
    &\Psi_3(a,1,b) \\&= 
    -(-1)^{|a|}f_1(a)f_2(1,b) + f_2(a,1)f_1(b) + 
    f_2(m_2(a,1),b) - f_2(a,m_2(1,b)) \\
    &= -0 + 0 + f_2(a,b)-f_2(a,b) = 0 \\
    &\Psi_3(a,b,1) \\&= 
    -(-1)^{|a|}f_1(a)f_2(b,1) + f_2(a,b)f_1(1) + 
    f_2(m_2(a,b),1) - f_2(a,m_2(b,1)) \\
    &= -0 + f_2(a,b) + 0 - f_2(a,b) = 0
  \end{align*}
  Hence, $\Psi_3=0$ whenever one input is a unit, and thus $m_3=0$
  when one input is a unit and we can choose $f_3=0$ when one input is
  a unit.

  Consider now some $n>3$. In the expression for $\Psi_n$, we have
  terms of the forms 
  \begin{align*} 
    &f_i(a_1,\dots,a_i)\cdot f_j(a_{i+1},\dots,a_n) & \text{and}\\
    &f_i(a_1,\dots,m_k(a_j,\dots,a_{j+k}),\dots,a_n) &.  
  \end{align*} 
  In the case that $a_1=1$ or $a_n=1$, the 1 occurs inside some $f_k$
  or $m_k$, with $k>1$, for all cases except for the terms
  $(-1)^{1+(n-2)|1|}f_1(1)f_{n-1}(a_2,\dots,a_n)$ and
  $(-1)^{0+2\cdot(\dots)}f_{n-1}(m_2(1,a_2),\dots,a_n)$. These have
  opposite signs, and their sum vanishes. Otherwise, in terms of the
  first kind, the unit must occur as an argument to one of the two
  $f_*$s, which consequently vanishes, and thus so does the entire
  term. Terms of the second kind vanish whenever the unit occurs
  outside the $m_k$. If the unit occurs within the $m_k$, then we
  distinguish between $k>2$ and $k=2$. If $k>2$, then $m_k$ vanishes
  by assumption. Thus, we only need to consider the case $k=2$. In
  this case two non-vanishing terms occur, namely
  \begin{multline*} 
    (-1)^{j} f_i(a_1,\dots,a_{i-1},m_2(1,a_{i+1}),\dots,a_n) = \\ 
    (-1)^{j} f_i(a_1,\dots,a_{i-1},a_{i+1},\dots,a_n) \\ 
  \end{multline*}
  \begin{multline*} 
    (-1)^{j-1} f_i(a_1,\dots,m_2(a_{i-1},1),a_{i+1},\dots,a_n) = \\ 
    -(-1)^{j} f_i(a_1,\dots,a_{i-1},a_{i+1},\dots,a_n)\quad.  
  \end{multline*}
  Hence these terms cancel each other and we conclude that
  $\Psi_n=0$. Thus $m_n=0$ follows and we can safely choose $f_n=0$. 

  Since the algorithm we derive does not use the internal structure of
  the proof of uniqueness of the operations, we refer the reader to
  \citep{kadeishvili80} for a proof. 
\end{proof}

This translates immediately into an algorithm for pointwise
computation of the $A_\infty$-structure maps.  For the computation of
an $A_\infty$-structure on $H_*A$, we need to fix some data for the
entire computation. Central to this is the choice of a cycle-choosing
map $f_1\colon H_*A\to A$. We need the map to send classes to cycles
representing the classes, but any such choice will work. This is not
the only choice made in the algorithm -- for each computed
$f_n(a_1,\dots,a_n)$, a choice is made of an appropriate coboundary
$dw$, to make $f_n$ a $\mathbb k$ module homomorphism. However, this
is easily seen to be reliable, as $H_*A$ is free over $\mathbb k$.

\begin{alg}[The Kadeishvili Algorithm]
  The algorithm takes as input a list of elements $a_1,\dots,a_n$ in
  $H_*A$, a cycle-choosing homomorphism $f_1\colon H_*A\to A$ as
  described above, and returns $m_n(a_1,\dots,a_n)$ and
  $f_n(a_1,\dots,a_n)$ fulfilling the Stasheff axioms $\St_n$ and
  $\St^m_n$.
  
  \begin{enumerate}
  \item If $n=1$, return $m_1(a_1)=0$ and $f_1(a_1)$ immediately.
  \item If $n=2$, set $\Psi_2(a_1,a_2)=f_1(a_1)f_1(a_2)$ and
    $m_2(a_1,a_2)=a_1a_2$ and go to step
    \ref{algo:nullhomotopy}. Otherwise, compute
    \begin{multline*}
      \Psi_n(a_1,\dots,a_n)= \\
      \sum_{s=1}^{n-1}(-1)^{\varepsilon_1(a_1,\dots,a_n,s)}
      f_s(a_1,\dots,a_s)\cdot f_{n-s}(a_{s+1},\dots,a_n) + \\
      \sum_{j=2}^{n-1}\sum_{k=0}^{n-j}(-1)^{\varepsilon_2(a_1,\dots,a_n,k,j)}
      f_{n-j+1}(a_1,\dots,a_k,m_j(a_{k+1},\dots,a_{k+j}),\dots,a_n)
    \end{multline*}
    where the expressions
    $\varepsilon_1(a_1,\dots,a_n,s)=s+(n-s+1)(|a_1|+\dots+|a_s|)$ and
    $\varepsilon_2(a_1,\dots,a_n,k,j)=k+j(n-k-j+|a_1|+\dots+|a_k|)$
    are the signs in the Stasheff morphism axiom $\St^m_n$ with the
    Koszul signs introduced.

    Note that the values of $f_k$ and $m_k$ for $k<n$ may be computed
    recursively using subsequent calls to this algorithm. The
    recursion bottoms out since $m_1$, $m_2$ and $f_1$ are already
    given.
  \item By the proof of the Minimality Theorem, the element
    $\Psi_n(a_1,\dots,a_n)\in A$ is a cycle. Hence, it belongs to some
    homology class $x$. Set $m_n(a_1,\dots,a_n)=x$.
  \item\label{algo:nullhomotopy} Since $m_n(a_1,\dots,a_n)$ is the
    homology class containing $\Psi_n(a_1,\dots,a_n)$, the
    cycle $f_1(m_n(a_1,\dots,a_n))$ is homologous to the
    cycle $\Psi_n(a_1,\dots,a_n)$. Thus
    $\Psi_n(a_1,\dots,a_n)-f_1(m_n(a_1,\dots,a_n))$ is a boundary, and
    we can pick an element $y$ such that
    $dy=\Psi_n(a_1,\dots,a_n)-f_1(m_n(a_1,\dots,a_n))$. We set
    $f_n(a_1,\dots,a_n)=y$, and return the higher multiplication
    $m_n(a_1,\dots,a_n)$ and the quasi-isomorphism component
    $f_n(a_1,\dots,a_n)$.
  \end{enumerate}
\end{alg}

\section{Computational reduction}
\label{sec:comp-reduct}

At a first glance, computing $A_\infty$-structures by a blackbox
method seems impractical due to the high degree of recursion and the
multiple infinities involved, i.e.\@ there are infinitely many arities
to compute, and even $H_*A$ tends to be infinite dimensional in group
cohomology, hence so is $(H_*A)^{\otimes n}$ for all the relevant
values of $n$.

In specific cases, however, we are able to reduce the complexity of
these computations to a managable size. For especially well-behaved
cohomology rings, we are able to reduce the computation of a full
$A_\infty$-structure to a finite problem.

We say that the induced $A_\infty$-algebra structures found by the
Minimality Theorem are $R$-linear if both the structure maps are
$R$-linear, and the quasi-isomorphism is an $R$-linear
$A_\infty$-morphism. 

\begin{thm}\label{lemma:periodic}
\theorema{lemma:periodic}
\end{thm}
\begin{proof}
  Set $\zeta=f_1(z)$. 
  We need to consider 
  \begin{multline*}
    \Psi_n(a_1,\dots,za_i,\dots,a_n) =
    \sum\pm f_j(a_1,\dots,za_i,\dots,a_j)f_{n-j}(a_{j+1},\dots,a_n) +
    \\
    \sum\pm f_j(a_1,\dots,a_j)f_{n-j}(a_{j+1},\dots,za_i,\dots,a_n) +
    \\
    \sum\pm
    f_{n-j+1}(a_1,\dots,za_i,\dots,m_j(a_{k+1},\dots,a_{k+j}),\dots,a_n) +
    \\
    \sum\pm
    f_{n-j+1}(a_1,\dots,m_j(a_{k+1},\dots,za_i,\dots,a_{k+j}),\dots,a_n) +
    \\
    \sum\pm
    f_{n-j+1}(a_1,\dots,m_j(a_{k+1},\dots,a_{k+j}),\dots,za_i,\dots,a_n)
    \quad .
    \\
  \end{multline*}
  In each summand of this expression, the term $za_i$ occurs within
  either a $f_j$ or a $m_j$ of lower arity than $n$. Hence, by
  assumption, we can commute $z$ out to a $\zeta$. Since, also,
  $\zeta$ commutes with all $f_n$ of lower arity, we find that
  \[
  \Psi_n(a_1,\dots,za_i,\dots,a_n) = \zeta\Psi_n(a_1,\dots,a_n)\quad.
  \]
  Hence $m_n(a_1,\dots,za_i,\dots,a_n)=zm_n(a_1,\dots,a_n)$
  follows. To finalize the argument, we need to find a boundary
  $h$ of the element $\zeta(\Psi_n-f_1m_n)(a_1,\dots,a_n)$ given a
  boundary $h'$ of $(\Psi_n-f_1m_n)(a_1,\dots,a_n)$. Note that
  \begin{align*}
    d(\zeta h') =
    (d\zeta) h' + \zeta dh' =
    \zeta(\Psi_n-f_1m_n)(a_1,\dots,a_n)\quad;
  \end{align*}
  since $\zeta=f_1(z)$, $f_1$ chooses cycles and hence
  $d\zeta=0$. Thus $\zeta h'$ is such a boundary.
\end{proof}

\begin{thm}\label{lemma:commutes}
  Suppose that $R$ is a finite $\mathbb{k}$-algebra and that 
  \begin{enumerate}[{B}1.]
  \item\label{lemma:commutes:1} $X$ is a periodic resolution of period
    $\pi$ of finitely generated $R$-modules, and that $A=\End_R(X)$ is
    the endomorphism dg-algebra of $X$.  Suppose further that there is
    some element $0\neq z\in H_*A$ such that we can choose
    $f_1(z)=\zeta$, a periodic map of period $\pi$ with each
    $\zeta_n=\Id$.
    
  \item\label{lemma:commutes:2} 
    for all $k<n$ we have constructed $m_k$ and $f_k$ such that
    A\ref{lemma:periodic:2} holds.
    
  \item\label{lemma:commutes:3} $b_1,\dots,b_t$ is a $\mathbb{k}[z]$-basis for
    $H_*A$ and for all $v_1,\dots,v_n$ chosen such that each
    $v_i\in\{b_1,\dots,b_t\}$ we know that $f_n(v_1,\dots,v_n)$ is
    periodic of period $\pi$ for all choices $v_1,\dots,v_n$.
  \end{enumerate}
  
  Then, from B\ref{lemma:commutes:1}, we can infer that $z$ generates
  a polynomial subalgebra of $H_*A$ and $H_*A$ is free over $\mathbb{k}[z]$.
  
  From B\ref{lemma:commutes:1}, B\ref{lemma:commutes:2} and the fact
  that $\zeta f_k(a_1,\dots,a_k)=f_k(a_1,\dots,a_k)\zeta$, 
  the homotopies $f_n(a_1,\dots,a_k)$ are periodic of period $\pi$.

  In addition, B\ref{lemma:commutes:3} implies that for all
  $a_1,\dots,a_n\in H_*A$, the map $f_n(a_1,\dots,a_n)$ is
  periodic of period $\pi$ and
  \begin{align*}
    m_n(a_1,\dots,za_i,\dots,a_n) &= zm_n(a_1,\dots,a_n) \\
    f_n(a_1,\dots,za_i,\dots,a_n) &= \zeta f_n(a_1,\dots,a_n) \\
    f_n(a_1,\dots,a_n)\zeta &= \zeta f_n(a_1,\dots,a_n)\quad.
  \end{align*}
\end{thm}
\begin{proof}
  If some $\zeta^N$ is null-homotopic, we can use the periodicity of
  $X$ to lower the degree of the null-homotopy, thereby inducing a
  null-homotopy for $\zeta$. However, we assumed that $z\neq0$. Hence
  $f_1(z)=\zeta$ is not null-homotopic. Thus, $\mathbb{k}[z]$ is a
  polynomial subalgebra of $H_*A$.  We set $I=H^{\geq 1}A$. This is an
  ideal in $H_*A$, and we can find the $\mathbb{k}$-vector space
  $J=I/I^2$ of indecomposables. We can pick a basis $b_1,\dots,b_r$ of
  $J/(z)$. Every $b_i$ has a representative in $J$, hence a
  representative that is not divisible by $z$. Furthermore,
  $b_1,\dots,b_r,z$ generate $H_*A$ as a $\mathbb{k}$-algebra. Suppose
  now that we had some dependency $\sum_i a_ib_i = 0$ over
  $\mathbb{k}[z]$. Then $z|a_i$ for all $a_i$, since otherwise the
  $b_i$ would not form a basis of $J/(z)$. But then we could divide
  the dependency by an appropriate power of $z$ and get a dependency
  involving the indecomposables. Hence $H_*A$ is free over
  $\mathbb{k}[z]$.

  From the condition $f_k(a_1,\dots,a_k)\zeta = \zeta
  f_k(a_1,\dots,a_k)$ we get by setting $d=|f_k(a_1,\dots,a_k)|$, that
  $(f_k(a_1,\dots,a_k)\zeta)_n=f_k(a_1,\dots,a_k)_n\zeta_{n+d}$ and
  that $(\zeta f_k(a_1,\dots,a_k))_n=\zeta_nf_k(a_1,\dots,a_k)_{n+\pi}$.
  Equality of chain maps forces the equality 
  $f_k(a_1,\dots,a_k)_n\zeta_{n+d}=\zeta_nf_k(a_1,\dots,a_k)_{n+\pi}$,
  and by the definition of $\zeta$, we are left with
  $f_k(a_1,\dots,a_k)_n=f_k(a_1,\dots,a_k)_{n+\pi}$.
  
  Since $b_1,\dots,b_t$ form a $\mathbb{k}[z]$-linear basis of $H_*A$, any
  element $a\in H_*A$ has a unique decomposition into a $\mathbb{k}[z]$-linear
  combination of the $b_i$.

  By Theorem \ref{lemma:periodic}, all the commutativity relations
  hold.

  For periodicity of $f_n(a_1,\dots,a_n)$, consider the terms of the
  difference $\Psi_n(a_1,\dots,a_n)-f_1m_n(a_1,\dots,a_n)$. Each term
  in this expression is either a composition of periodic maps of
  period $\pi$, or a periodic map of period $\pi$, by the assumptions on
  all $f_k$. Hence, $\Psi_n(a_1,\dots,a_n)-f_1m_n(a_1,\dots,a_n)$ is
  periodic of period $\pi$. 

  Finally, by assumption B\ref{lemma:commutes:3}, $f_n(v_1,\dots,v_n)$
  is periodic of period $\pi$, for all choices of
  $v_1,\dots,v_n\in\{b_1,\dots,b_t\}$. Hence
  $f_n(v_1,\dots,zv_k,\dots,v_n)=\zeta f_n(v_1,\dots,v_n)$, and by
  Theorem \ref{lemma:periodic}, the resulting homotopy
  $f_n(a_1,\dots,a_n)$ is given by composing some $\zeta^s$, which has
  period $\pi$, with $f_n(v_1,\dots,v_n)$, which is also periodic of
  period $\pi$.
\end{proof}

\begin{thm}\label{lemma:halts}
\theoremb
\end{thm}
\begin{proof}
  The proof follows by induction. Suppose that $\kappa > 2q-2$, and that
  we have already proven $f_k=0$ and $m_k=0$ for all $q\leq
  k<\kappa$. In the computational step where we compute $f_\kappa$ and
  $m_\kappa$, we start by considering $\Psi_\kappa$. This expression
  has two kinds of terms.

  First, there are the terms of the form $f_i\cdot
  f_{\kappa-i}$. Since $\kappa>2q-2$, either $i\geq q$ or $\kappa-i\geq
  q$. Hence by the induction hypothesis, $f_i\cdot f_{\kappa-i}=0$.
  Second, there are the terms of the form $f_i\circ_j
  m_{\kappa-i+1}$. Again, either $\kappa-i+1\geq\kappa-i\geq q$ or
  $i\geq q$. Hence, by hypothesis either $f_i=0$ or
  $m_{\kappa-i+1}=0$.
  Hence $\Psi_\kappa=0$. Thus we can choose $m_\kappa=0$ and
  $f_\kappa=0$. This proves the induction step and completes the
  proof.
\end{proof}

\subsection{Application: Cohomology of cyclic groups}
\label{sec:appl-cohom-cycl}

Suppose $\mathbb{k}$ is a finite field and
$R=\mathbb{k}[\alpha]/(\alpha^q)$. Then
$\Ext_R^*(\mathbb{k},\mathbb{k})$ is, for $q$ a power of the
characteristic $p$ of $\mathbb{k}$, the ring
$H^*(C_q,\mathbb{k})$. For $q\geq 4$, it has a finite presentation, as
a ring, given by $\Lambda^*(x)\otimes \mathbb{k}[y]$. We can choose a
minimal periodic free resolution $(P_*,d)$ of $\mathbb{k}$ using
$R$-modules on the form
\begin{diagram}
\cdots &\rTo& R 
&\rTo^{\cdot \alpha}& R
&\rTo^{\cdot \alpha^{q-1}}& R
&\rTo^{\cdot \alpha}& R
&\rTo^{\cdot \alpha^{q-1}}& R
&\rTo^{\cdot \alpha}& R
&\rTo& R \to \mathbb{k} \to 0
\end{diagram}

Now, $\Ext_R^*(\mathbb{k},\mathbb{k})$ is the homology of the
dg-algebra of graded module maps $P_*\to P_*$, with the induced
differential $\partial f=df-(-1)^{|f|}fd$, and thus we can find
representatives for the classes $x$ and $y$ given by $\xi$ and $\eta$:

\begin{diagram}
  \dots &\rTo& 
  \mathbb{k}G &\rTo^{\cdot \alpha^{q-1}}&
  \mathbb{k}G &\rTo^{\cdot\alpha}&
  \mathbb{k}G &\rTo^{\cdot \alpha^{q-1}}&
  \mathbb{k}G &\rTo^{\cdot\alpha} & \mathbb{k}G \to \mathbb{k} \to 0\\
  \xi: &&
  \dTo_{\cdot\alpha^{q-2}} &&    
  \dTo_{\cdot 1} &&    
  \dTo_{\cdot\alpha^{q-2}} &&
  \dTo_{\cdot 1} & \rdTo^\epsilon&    &.
  \\
  \dots &\rTo& 
  \mathbb{k}G &\rTo^{\cdot\alpha}&
  \mathbb{k}G &\rTo^{\cdot \alpha^{q-1}}&
  \mathbb{k}G &\rTo^{\cdot\alpha}&
  \mathbb{k}G &\rTo^\epsilon & \mathbb{k} &\rTo& 0 \\
\end{diagram}

\begin{diagram}
  \dots &\rTo& 
  \mathbb{k}G &\rTo^{\cdot\alpha}&
  \mathbb{k}G &\rTo^{\cdot \alpha^{q-1}}&
  \mathbb{k}G &\rTo^{\cdot\alpha}& 
  \mathbb{k}G &\rTo^{\cdot \alpha^{q-1}}& 
  \mathbb{k}G &\rTo^{\cdot\alpha}& \mathbb{k}G& \to \mathbb{k} \to 0 \\
  \eta: &&
  \dTo_{\cdot 1} &&    
  \dTo_{\cdot 1} &&    
  \dTo_{\cdot 1} &&
  \dTo_{\cdot 1} & \rdTo^\epsilon&   &&&&&. \\
  \dots &\rTo& 
  \mathbb{k}G &\rTo^{\cdot\alpha}&
  \mathbb{k}G &\rTo^{\cdot \alpha^{q-1}}&
  \mathbb{k}G &\rTo^{\cdot\alpha}&
  \mathbb{k}G &\rTo^\epsilon & \mathbb{k} &\rTo& 0\\
\end{diagram}

The map $\eta$ fulfills all the requirements for $\zeta$ in Theorem
\ref{lemma:commutes}, and hence, as long as all computed homotopies
have periodicity dividing 2, it will suffice to consider parameter
sets taken from $\{1,x\}$ when computing the higher
multiplications. Furthermore, since we may compute a strictly unital
$A_\infty$-algebra structure, the only parameter sets we need to
consider are of the form $x\otimes\dots\otimes x$. By
$\mathbb{k}[y]$-linearity, all other values follow. The
$A_\infty$-algebra structure on this cohomology ring computed by Dag
Madsen \citep{madsen_phd} agrees with the $A_\infty$-algebra structure
found using the techniques in this paper.

Performing the Kadeishvili algorithm in this setting yields the maps
$f_m(x,\dots,x)$ as given by $h$ in:
\begin{diagram}
  \dots &\rTo& 
  \mathbb{k}G &\rTo^{\cdot \alpha^{q-1}}&
  \mathbb{k}G &\rTo^{\cdot\alpha}&
  \mathbb{k}G &\rTo^{\cdot \alpha^{q-1}}&
  \mathbb{k}G &\rTo^{\cdot\alpha} & \mathbb{k}G \\
  h: &&
  \dTo_{\cdot(-\alpha^{q-1-m})} &&    
  \dTo_{\cdot 0} &&    
  \dTo_{\cdot(-\alpha^{q-1-m})} &&
  \dTo_{\cdot 0} & \rdTo^\epsilon&    &.
  \\
  \dots &\rTo& 
  \mathbb{k}G &\rTo^{\cdot\alpha}&
  \mathbb{k}G &\rTo^{\cdot \alpha^{q-1}}&
  \mathbb{k}G &\rTo^{\cdot\alpha}&
  \mathbb{k}G &\rTo^\epsilon & \mathbb{k} \\
\end{diagram}
and correspondingly $m_m(x,\dots,x)=0$ for all $m<q-1$. Then
$f_{q-1}(x,\dots,x)$ is given by $h$ in:
\begin{diagram}
  \dots &\rTo& 
  \mathbb{k}G &\rTo^{\cdot \alpha^{q-1}}&
  \mathbb{k}G &\rTo^{\cdot\alpha}&
  \mathbb{k}G &\rTo^{\cdot \alpha^{q-1}}&
  \mathbb{k}G &\rTo^{\cdot\alpha} & \mathbb{k}G \\
  h: &&
  \dTo_{\cdot(-1}) &&    
  \dTo_{\cdot 0} &&    
  \dTo_{\cdot(-1)} &&
  \dTo_{\cdot 0} & \rdTo^\epsilon&    &.
  \\
  \dots &\rTo& 
  \mathbb{k}G &\rTo^{\cdot\alpha}&
  \mathbb{k}G &\rTo^{\cdot \alpha^{q-1}}&
  \mathbb{k}G &\rTo^{\cdot\alpha}&
  \mathbb{k}G &\rTo^\epsilon & \mathbb{k} \\
\end{diagram}
and thus $m_q(x,\dots,x)=y$ and $f_q(x,\dots,x)=0$. Computing further
will reveal $f_m(x,\dots,x)=0$ and $m_m(x,\dots,x)=0$ for all $q+1\leq
m\leq 2q$, which completes the computation of an
$A_\infty$-algebra structure on $\Ext_R^*(\mathbb{k},\mathbb{k})$.

\section{Implementation}
\label{sec:implementation}

This approach is implemented in the \textsc{Magma} computer algebra
system \citep{magma} as a component in the computational group
cohomology modules. The implementation of the Kadeishvili algorithm
expects to be working with some $\End_R(X)$, but the supporting
homological algebra functionality was built to compute group
cohomology rings. Hence, currently, the module will only work smoothly
computing $A_\infty$-algebra structures on cohomology rings of finite
groups.

The user interface has three functions at its core:
\begin{description}
\item[\texttt{AInfinityRecord(G,n)}] Yields a computation object
  storing internal information for the computations. Expects a finite
  group $G$ and the length of a partial projective resolution of the
  trivial module $\mathbb{k}$ over the group algebra $\mathbb{k}G$.
\item[\texttt{HighProduct(Aoo,lst)}] Computes
  $m_k(a_1,\dots,a_k)$. Expects a computation object \texttt{Aoo}, as
  produced by \texttt{AInfinityRecord} and a sequence \texttt{lst} of
  length $k$ of elements of the cohomology ring stored as
  \texttt{Aoo`S}. Returns $m_k(a_1,\dots,a_k)$ as an element of
  \texttt{Aoo`S}.
\item[\texttt{HighMap(Aoo,lst)}] Computes the quasi-isomorphism
  component $f_k(a_1,\dots,a_k)$. Expects a computation object
  \texttt{Aoo}, as produced by \texttt{AInfinityRecord} and a sequence
  \texttt{lst} of length $k$ of elements of \texttt{Aoo`S}. Returns
  $f_k(a_1,\dots,a_k)$ as a chain map endomorphism of the partial
  projective resolution stored in \texttt{Aoo`P}.
\end{description}

\appendix

\bibliographystyle{abbrvnat}
\bibliography{library}

\end{document}